\documentclass[a4paper,12pt]{article}
\usepackage[T2A]{fontenc}
\usepackage[utf8]{inputenc}
\usepackage{amssymb}
\usepackage{amsmath}

\newcounter{theorem}
\newenvironment{theorem}{ %%
	\refstepcounter{theorem}%%
	\par\medskip\noindent{\bfseries Theorem \thetheorem.} %%
	\itshape %%
}{ %%
	\upshape\par %%
}
\renewcommand{\ge}{\geqslant}
\renewcommand{\le}{\leqslant}

\newcommand{\Z}{\mathbb{Z}}
\newcommand{\N}{\mathbb{N}}
\newcommand{\nd}{\nobreakdash-} 
\begin{document}

\author{Alexey~Bashtanov}
\title{Generic mixing is rank\nd $1$\footnote{Supported by Leading Scientific School Support Grant NSh-5998.2012.1}}
\date{29 Mar 2012}

\maketitle

\begin{abstract}
In paper~\cite{tikhonov0} S.~V.~Tikhonov introduced a complete metric on the space of mixing transformations. It generates a topology called leash-topology.
In~\cite{tikhonov1} Tikhonov states the following problem: for what mixing transformation $T$ its conjugacy class is dense in the space of mixing transformations 
with leash-topology?
We show that it is true for every mixing $T$. As a corollary we get that generic mixing is rank\nd $1$.
\end{abstract}

\section{Conjugacy classes density}
\paragraph{Preliminaries.}
Invertible measure-preserving transformation $T$ acting on a standard non-atomic probability space\footnote{Standard non-atomic probability space is a probability 
space isomorphic $mod$ $0$ to unit interval with Lebesgue measure. We consider transformations and other relations equal if they differ on zero measure set.} 
$(X, \Sigma, \mu)$ is called {\em mixing} if for any measurable sets $A$, $B$
$$\lim_{n\to\infty}\mu(T^nA\cap B)=\mu(A)\mu(B).$$

Let $\{A_i\}$ be a countable collection of sets generating $\sigma$\nd algebra $\Sigma$.
Each one of the metrics
$$d(T,S) = \sum_{i\in\N}\frac 1 {2^i}(\mu(TA_i\triangle SA_i) + \mu(T^{-1}A_i\triangle S^{-1}A_i))$$
and 
$$a(T,S) = \sum_{i,j\in\N}\frac 1 {2^{i+j}}|\mu(TA_i\cap A_j) - \mu(SA_i\cap A_j) |$$
generates weak topology in the space~$\mathcal M$ of mixing transformations. 

None of these metrics is complete in~$\mathcal M$. To make~$\mathcal M$ a complete metric space 
Tikhonov introduced the metric
$$\tau(T,S) = d(T,S) + \sup_{i\in\Z} a(T^n, S^n)$$
in his paper~\cite{tikhonov0}. This metric generates a different topology called {\em leash-topology}. The collection of neighborhoods of form
$$\mathcal U(T, q, \varepsilon) = \{S\in\mathcal M\mid\forall n\in\Z\;\forall A,B\in q\;|\mu(T^nA\cap B) - \mu(S^nA\cap B)|<\varepsilon\}$$
form a base of this topology. Each neighborhood of this form is defined by its center $T$, a collection $q$ of measurable sets and a positive number $\varepsilon$.

One can explore generic properties of mixings in the space~$\mathcal M$, equipped with the metric $\tau(T,S)$.
A subset of metric space is called a {\em $G_\delta$ set} if it can be represented as a countable intersection of open sets.
A superset of a dense $G_\delta$ set is called {\em residual}. If a set of transformations having some property is residual, the property is called {\em generic}.
Term ``generic'' explains itself with the following facts: countable intersection of residual sets is residual, a superset of a residual set is residual.

The main result of this paper is
\begin{theorem}\label{main}
Conjugacy class of any mixing transformation $S$ is dense in~$\mathcal M$.
\end{theorem}

\paragraph{Bernoulli shift with entropy $1$.}
Let us define a mixing called {\em Bernoulli shift with entropy\footnote{For notational convenience we use binary entropy function. If one uses natural logarithm, entropy of
this transformation would be $\ln 2$}~$1$}. 
Its phase space $X$ is a Cartesian product of countably many copies of space $\{0, 1\}$:
 $$X=\{0,1\}^\Z, \mu(\{0\})=\mu(\{1\})=\frac 1 2.$$
$X$ is a standard non-atomic probability space. Its elements are bidirectionally infinite zero-one sequences, and the transformation acts on them by right shift.
We will write points of the space as $$x=(\dots x_{-2}, x_{-1}, x_0, x_1, x_2, \dots).$$
Let us denote the sets $\{x\in X\mid x_0 = 0\}$ and $\{x\in X\mid x_0 = 1\}$ as $D_0$ and $D_1$ respectively. Let us denote the partition $X=D_0\sqcup D_1$ as $\xi$.
This partition is generating: the intersection $$\bigvee_{i\in\Z}T^i\xi$$ of its images is the partition into points.
This implies every measurable set $C\subset X$ can be approximated by sets $\{C_k\}_{k\in\N}$, where $C_k$ is measurable with respect to $$\bigvee_{|i|\le k}T^i\xi.$$

Tikhonov showed~(\cite{tikhonov1}) that Bernoulli shifts with entropy~$1$ are dense in~$\mathcal M$.
It follows from a more general statement by Tikhonov: transformations conjugate to any fixed Cartesian product are dense in $\mathcal M$.
Bernoulli shift is always a Cartesian product of another Bernoulli shifts with lesser entropy.

Thus, to prove that conjugacy class of any mixing transformation is dense in $\mathcal M$ it is enough to show that 
in every neighborhood of any Bernoulli shift with entropy $1$ there exists an element of this conjugacy class.

\paragraph{Passing to smaller neighborhoods.}
Let us fix $T$~--- a Bernoulli shift with entropy~$1$, and its neghtbourhood $\mathcal U(T, q, \varepsilon)$. 
Let us denote the elements of the collection $q$ by $A_1, A_2, \dots A_n$.
For every $A_j$ let $\tilde A_j$ be measurable with respect to some finite intersection 
$$\bigvee_{|i|\le k_j}T^i\xi$$
and let $\tilde A_j$ approximate $A_j$ with precision $\frac\varepsilon 5$.

Denote the collection $\{\tilde A_1, \tilde A_2, \dots \tilde A_n\}$ as $\tilde q$, and the maximum of $k_j$ as $k$.
Every~$\tilde A_i$ is measurable with respect to $$\bigvee_{|i|\le k}T^i\xi,$$ and for every $j$ the sets $A_j$ and $\tilde A_j$ are close to each other: 
$\mu(\tilde A_j\triangle A_j)<\frac\varepsilon 5$.

Let us show that $$\mathcal U(T, \tilde q, \frac\varepsilon 5)\subset\mathcal U(T, q, \varepsilon).$$
Indeed, if $P\in\mathcal U(T, \tilde q, \frac\varepsilon 5)$, then for any $r, s$ and $m$
$$|\mu(P^m\tilde A_r\cap\tilde A_s) - \mu(T^p\tilde A_r\cap \tilde A_s)|<\frac\varepsilon 5.$$
This implies that for any $r, s$ and $m$
\begin{multline*}
|\mu(P^mA_r\cap A_s) - \mu(T^mA_r\cap A_s)|\le\\\le
  |\mu(P^mA_r\cap A_s) - \mu(P^m\tilde A_r\cap A_s)| + |\mu(P^m\tilde A_r\cap A_s) - \mu(P^m\tilde A_r\cap\tilde A_s)| +\\+
  |\mu(P^m\tilde A_r\cap\tilde A_s) - \mu(T^m\tilde A_r\cap \tilde A_s)| +\\+
  |\mu(T^m\tilde A_r\cap\tilde A_s) - \mu(T^m\tilde A_r\cap A_s)| + |\mu(T^m\tilde A_r\cap A_s) - \mu(T^mA_r\cap A_s)|\le\\\le
  \frac\varepsilon 5 + \frac\varepsilon 5 + \frac\varepsilon 5 + \frac\varepsilon 5 + \frac\varepsilon 5 =\varepsilon,
\end{multline*}
that is $P\in\mathcal U(T, q, \varepsilon)$.

Let us consider now a collection $\hat q$, consisting of all atoms of partition $$\bigvee_{|i|\le k}T^i\xi.$$
Every atom is of the form $$B_{(b_{-k}, b_{-k+1}, \dots b_k)} = \bigcap_{|i|\le k} T^iD_{b_i} =
 \{x\in X\mid(x_{-k}=b_{-k}, x_{-k+1}=b_{-k+1}, \dots x_k=b_k\}$$

Let us show that $$\mathcal U(T, \hat q, \frac\varepsilon{5\cdot 2^{4k+2}})\subset\mathcal U(T, \tilde q, \frac\varepsilon 5).$$
Indeed, let $P\in\mathcal U(T, \hat q, \frac\varepsilon{5\cdot 2^{4k+2}})$, that is for all $u, v$ and $m$
$$|\mu(P^mB_u\cap B_v) - \mu(T^mB_u\cap B_v)|<\frac\varepsilon{5\cdot 2^{4k+2}}.$$
Here $u$ and $v$ are any one-zero sequences of length $2k+1$.
To estimate the difference $|\mu(P^m\tilde A_r\cap\tilde A_s) - \mu(T^m\tilde A_r\cap\tilde A_s)|$ notice that $\tilde A_r$ 
can be pieced out of some $C_{u_1}, C_{u_2}, \ldots\in\hat q$, and the number of the pieces cannot exceed $2^{2k+1}$.
Similarly decompose $\tilde A_s$ into pieces $C_{v_1}, C_{v_2}, \ldots\in\hat q$.
\begin{multline*}
|\mu(P^m\tilde A_r\cap\tilde A_s) - \mu(T^m\tilde A_r\cap\tilde A_s)| = \Bigl|\sum_{l,w}\mu(P^mC_{u_l}\cap C_{v_w}) - \sum_{l,w}\mu(T^mC_{u_l}\cap C_{v_w})\Bigr|\le\\\le
\sum_{l,w}|\mu(P^mC_{u_l}\cap C_{v_w}) - \mu(T^mC_{u_l}\cap C_{v_w})|<\sum_{l,w}\frac\varepsilon{5\cdot 2^{4k+2}}\le\\\le
2^{2k+1}\cdot 2^{2k+1}\cdot \frac\varepsilon{5\cdot 2^{4k+2}}=\frac\varepsilon 5.
\end{multline*}

\paragraph{Images almost-independence.}
Recall $S$ is an arbitrary mixing. 
Let $S$ act on standard non-atomic probability space~$Y$. Let us show that in the neigbourhood $\mathcal U(T, \hat q, \frac\varepsilon{5\cdot 2^{4k+2}})$
there is a transformation $V$ acting on $X$ and conjugate to $S$. According to~\cite{bashtanov} for every $\delta>0$ there exists such a set $A\subset Y$ of measure~$\frac 1 2$ that
sets $T^mA, m\in\Z$ are collectionwise $\delta$\nd independent. Condition of collectionwise {\em $\delta$\nd independence} of sets $E_1, E_2, \dots$ means that
the measure of any finite intersection of form
\begin{equation}
\label{indep}
E_{m_1}\cap E_{m_2}\cap\dots E_{m_p}
\end{equation}
 differs from $$\prod_{i=1}^p\mu(E_{m_i})$$ less than by $\delta$ in assumption that $m_i$ are different. 
To build such a set in~\cite{bashtanov} author uses the method offered by V.~V.~Ryzhikov in~\cite{ryzhikovPairwise}:
one builds an approximating sequence of Rokhlin castles, random union of cells of sufficiently big castle satisfies the conditions with positive probability. 

Let us call sets $E_1, E_2, \ldots \subset Y$ {\em well collectionwise $\delta$\nd independent}, if
the measure of any finite intersection of form
\begin{equation}
\label{full_indep}
E_{m_1}\cap E_{m_2}\cap\dots E_{m_{p_0}}\cap(Y\setminus E_{n_1})\cap(Y\setminus E_{n_2})\cap\dots(Y\setminus E_{n_{p_1}})
\end{equation}
differs from $$\prod_{i=1}^{p_0}\mu(E_{m_i})\prod_{i=1}^{p_1}\mu(Y\setminus E_{n_i})$$ less than by $\delta$ in assumption that none of $m_i$ and $n_i$ coincide. It is not hard to show that if a collection of sets
is $\delta$\nd independent, than every its finite subcollection of cardinality $c$ will be well collectionwise $c\delta$\nd independent. This follows from the fact that
any intersection of form \eqref{full_indep} and length~$c$ can be expressed by 
$c$ intersections of form~\eqref{indep} using operations of disjoint union and set-theoretical difference.
Notice that using the methods of~\cite{ryzhikovPairwise, bashtanov} one can directly prove well collectionwise $\delta$\nd independence of all images of the set $A$.

\paragraph{Building the conjugate transformation.}
Let us take $\delta = \frac 1{4k+2}\cdot\frac\varepsilon{15\cdot 2^{4k+2}}$. 
Then every $4k+2$ images of set $A$ will be well $\frac\varepsilon{15\cdot 2^{4k+2}}$\nd independent.
Denote $A$ as $F_0$, $X\setminus A$ as $F_1$, and the partition $Y = F_0\sqcup F_1$ as $\eta$. Consider atoms of partitions $$\eta_k = \bigvee_{|i|\le k}T^i\eta.$$
The measure of each of them differs from $\frac 1{2^{2k+1}}$ less than by $\frac\varepsilon{15\cdot 2^{4k+2}}$.
In $X$ there exists a partition $$\xi_k = \bigvee_{|i|\le k}T^i\xi,$$ and the measure of each its atom is exactly $\frac 1{2^{2k+1}}$.
The desired transformation is $V = Q^{-1}SQ$, where $Q\colon X\to Y$ is an arbitrary measure-preserving invertible transformation, that maps atoms of partition
$\xi_k$ onto the correspondent atoms of $\eta_k$ with a gap less than $\frac\varepsilon{15\cdot 2^{4k+2}}$ for each atom. 

\paragraph{Estimates for the transformation built.}
Let us prove that the transformation $V$ is situated in the neighborhood $\mathcal U(T, \hat q, \frac\varepsilon{5\cdot 2^{4k+2}})$.
For this it is necessary for every $u, v$ and $m$ the following condition to satisfy:
$$|\mu(V^mB_u\cap B_v) - \mu(T^mB_u\cap B_v)|<\frac\varepsilon{5\cdot 2^{4k+2}}.$$
Let us represent $B_u$ and $B_v$ in such a form:
$$B_u = \bigcap_{|i|\le k} T^iD_{u_i},\quad B_v = \bigcap_{|i|\le k} T^iD_{v_i}.$$
Then $T^mB_u\cap B_v$ and $V^mB_u\cap B_v$ will take form
$$\bigl(\bigcap_{|i|\le k} T^{i+m}D_{u_i}\bigr)\cap\bigl(\bigcap_{|i|\le k} T^iD_{v_i}\bigr)
\quad\mbox{and}\quad
\bigl(\bigcap_{|i|\le k} V^mT^iD_{u_i}\bigr)\cap\bigl(\bigcap_{|i|\le k} T^iD_{v_i}\bigl)$$
respectively. Estimate the difference of measures: 
\begin{multline*}
\biggl|\mu\Bigl(\bigl(\bigcap_{|i|\le k} T^{i+m}D_{u_i}\bigr)\cap\bigl(\bigcap_{|i|\le k} T^iD_{v_i}\bigr)\Bigr) - 
      \mu\Bigl(\bigl(\bigcap_{|i|\le k} V^mT^iD_{u_i}\bigr)\cap\bigl(\bigcap_{|i|\le k} T^iD_{v_i}\bigr)\Bigr)\biggr| \le\\\le
\biggl|\mu\Bigl(\bigl(\bigcap_{|i|\le k} T^{i+m}D_{u_i}\bigr)\cap\bigl(\bigcap_{|i|\le k} T^iD_{v_i}\bigr)\Bigr) - 
      \mu\Bigl(\bigl(\bigcap_{|i|\le k} S^{i+m}F_{u_i}\bigr)\cap\bigl(\bigcap_{|i|\le k} S^iF_{v_i}\bigr)\Bigr)\biggr| +\\+
\biggl|\mu\Bigl(\bigl(\bigcap_{|i|\le k} S^{i+m}F_{u_i}\bigr)\cap\bigl(\bigcap_{|i|\le k} S^iF_{v_i}\bigr)\Bigr) -
       \mu\Bigl(V^m\bigl(\bigcap_{|i|\le k} T^iD_{u_i}\bigr)\cap\bigl(\bigcap_{|i|\le k} T^iD_{v_i}\bigl)\Bigl)\biggl|.
\end{multline*}
Let us estimate the first item. Its left side is some intersection of $D_0$, $D_1$ and their images under action by powers of $T$, not more than
 $4k+2$ sets are intersected. 
They are collectionwise independent, or, what is the same, well collectionwise $0$\nd independent.
The right side is the similar intersection of $F_0$, $F_1$ and their images under action by powers of $S$. They are well 
collectionwise $\frac\varepsilon{15\cdot 2^{4k+2}}$\nd independent.
Hence the difference between the left and the right sides does not exceed $\frac\varepsilon{15\cdot 2^{4k+2}}$ in modulus.
It can happen that some sets of the left side coincide. In this case the correspondent sets on the right also coincide and we have the intersection of less number of sets.
Also there can be an intersection of form $$\cdots\cap T^zD_0\cap\cdots\cap T^zD_1\cap\cdots$$ on the left. Then there is a similar intersection 
$$\cdots\cap S^zF_0\cap\cdots\cap S^zF_1\cap\cdots$$ 
on the right. In this case both left and right sides equal $0$.

Let us estimate the second item. 
\begin{multline*}
\biggl|\mu\Bigl(\bigl(\bigcap_{|i|\le k} S^{i+m}F_{u_i}\bigr)\cap\bigl(\bigcap_{|i|\le k} S^iF_{v_i}\bigr)\Bigr) -
       \mu\Bigl(V^m\bigl(\bigcap_{|i|\le k} T^iD_{u_i}\bigr)\cap\bigl(\bigcap_{|i|\le k} T^iD_{v_i}\bigl)\Bigl)\biggl|=\\=
\biggl|\mu\Bigl(S^m\bigl(\bigcap_{|i|\le k} S^iF_{u_i}\bigr)\cap\bigl(\bigcap_{|i|\le k} S^iF_{v_i}\bigr)\Bigr) -
       \mu\Bigl(S^mQ\bigl(\bigcap_{|i|\le k} T^iD_{u_i}\bigr)\cap Q\bigl(\bigcap_{|i|\le k} T^iD_{v_i}\bigl)\Bigl)\biggl|
\end{multline*}
Each of the sets $$\bigcap_{|i|\le k} S^iF_{u_i}\mbox{ and }\bigcap_{|i|\le k} S^iF_{v_i}$$ on the left differs from the corresponding set
$$Q\bigcap_{|i|\le k}T^iD_{u_i}\mbox{ or }Q\bigcap_{|i|\le k}T^iD_{v_i}$$ on the right by a set of measure less than ~$\frac\varepsilon{15\cdot 2^{4k+2}}$.
Therefore the difference between the left and the right sides is not more than $\frac{2\varepsilon}{15\cdot 2^{4k+2}}$.

Sum of the items is less than $\frac\varepsilon{15\cdot 2^{4k+2}} + \frac{2\varepsilon}{15\cdot 2^{4k+2}} = \frac\varepsilon{5\cdot 2^{4k+2}}$, as was to be proved. 

We have built a transformation $V$, conjugate to $S$ and situated in the neighborhood $\mathcal U(T, \hat q, \frac\varepsilon{5\cdot 2^{4k+2}})$, and therefore also in $\mathcal U(T, q, \varepsilon)$.
The theorem is proved.

\section{Genericity of rank\nd $1$}
\paragraph{Preliminaries.}
Let the transformation $T$ act on the space $(X, \Sigma, \mu)$. A partition $\xi$ of form
$$X = E\sqcup TE\sqcup\dots\sqcup T^{n-1}E\sqcup D$$
is called a {\em Rokhlin tower of height $n$}.
Sets $E, TE, \dots T^{n-1}E$ are called levels of the tower, and $D$ is called a remainder.
Let us say that {\em a tower approximates} set $A$ with accuracy $\varepsilon > 0$, if there exists such a union $B$ of its levels, 
that $\mu(A\triangle B)<\varepsilon$. Let us say that {\em a sequence of towers approximates} set $A$ with accuracy $\varepsilon > 0$,
if all the towers in the sequence starting from a certain index approximate the set with the specified accuracy.
If there exists such a sequence of towers $\{\xi_i\}_{i\in\N}$ that 
for every measurable set $A$ and number $\varepsilon>0$ the sequence $\{\xi_i\}$ approximates $A$ with accuracy $\varepsilon$, then
$T$ is said to be {rank\nd $1$}. Rank\nd $1$ mixing transformations were first introduced by D.~Ornstein in~\cite{ornstein}.

\paragraph{Residuality of the set.}
Let us prove the genericity of rank\nd $1$ in mixings.
\begin{theorem}
The set of rank\nd $1$ mixings is a dense $G_\delta$ subset in $\mathcal M$.
\end{theorem}
Let us fix a countable collection of sets $\{A_j\}$, such that every measurable set $A\subset X$ can be approximated by elements of the collection
 (such a collection exists due to the properties of standard probability space).
If a sequence of towers approximate $A_j$ with arbitrary little accuracy, then the transformation is rank\nd $1$.
Let us denote the set of mixings, having a Rokhlin tower approximating the sets $A_1, A_2, \dots A_j$ with accuracy $a$, where $a<\frac 1 k$, as $\mathcal R(j, k)$.
Then the intersection $$\bigcap_{j,k}\mathcal R(j, k)$$
is exactly the set of rank\nd $1$ mixing. 

Let us show that every $\mathcal R(j, k)$ is open. Let $S\in\mathcal R(j, k)$. It means that there exists such a tower $\xi$,
that for every $i\in\{1, 2, \dots j\}$ the set $A_i$ differs from some union $B_i$ of tower levels by a set of measure $a<\frac 1 k$. 
Let $E$ be the downmost level of tower~$\xi$ and $n$ be its height.
Denote $$b=\frac 1{n^2}\left(\frac 1 k - a\right).$$
Let a mixing $V$ be situated in a neighborhood $\mathcal U(S, \{E, SE, \dots S^{n-1}E\}, b)$. 

Consider the set $$\tilde E = E\setminus(VE\cup V^2E\cup\dots\cup V^{n-1}E).$$ It differs from $E$ not more than by measure $(n-2)b$.
By construction of $\tilde E$ it is a downmost level of some tower $\eta$: $$X=\tilde E\sqcup V\tilde E\sqcup\dots\sqcup V^{n-1}\tilde E\sqcup \tilde D.$$
As far as the downmost level of the tower $\eta$ differs from the downmost level of tower $\xi$ not more than by measure $(n-2)b$, the other corresponding levels
of these two towers differ not more than by measure $b+(n-2)b = (n-1)b$ one from another, and the corresponding level unions differ not more than by measure $n(n-1)b$.
A union of levels of tower $\eta$ corresponding to some level union $B_i$ in $\xi$ differs from $B_i$ not more than by measure $n(n-1)b$, so
it differs from $A_i$ not more than by measure $n(n-1)b+a < \frac 1 k$. 

We achieved that an arbitrary transformation $S$ is included in $\mathcal R(j, k)$ with some its neighborhood. That is, $\mathcal R(j, k)$ is open, 
therefore the set of rank\nd $1$ mixings is a $G_\delta$ set.
Conjugacy class of any rank\nd $1$ mixing is dense in $\mathcal M$, so the set of rank\nd $1$ mixing is residual.
\paragraph{Corollaries.}
As far as rank\nd $1$ is generic for mixing, all the properties of rank\nd $1$ mixing are generic. 
Expressly minimal self-joinings of all orders, and all its corollaries: trivial centralizer, primeness and absence of roots (see~\cite{king}).

\section{Remarks}
\paragraph{Approximating by conjugates of a non-mixing transformation.}
In paper~\cite{bashtanov} author uses a more weak property than mixing~--- absence of partial rigidity.
Transformation $T$ is called {\em partially rigid with coefficient $a>0$}, if for every measurable set $A\subset X$
$$\limsup_i\mu(A\cap T^iA) \ge a\mu(A).$$
Let us expand the space $\mathcal M$ of mixings to the space $\mathcal A$ of all invertible measure-preserving transformations.
We will use the same metric $\tau(T, S)$. 
It is interesting that in this space there exists such an analog to theorem~\ref{main}:
\begin{theorem}
Fix a transformation $S$, not partially rigid. Then in any neighborhood of any mixing $T\in\mathcal A$ there exists a 
transformation, conjugate to~$S$.
\end{theorem}
Proof is analogous. Note that the opposite is true: conjugacy class of transformation $S$ cannot approximate mixings if $S$ is partially rigid.
\paragraph{Direct proof of density of rank one.}
Ryzhikov informed author that he knows another proof that rank\nd $1$ mixings are dense in $\mathcal M$. This proof is based neither on Tikhonov theorem about density of 
conjugates to a Cartesian product nor on results \cite{ryzhikovPairwise, bashtanov}.

\end{document}